\documentclass[a4paper,10pt]{elsarticle}

\usepackage{amsmath, amssymb, latexsym,enumerate, mathrsfs}

\usepackage{lineno,hyperref}

\newtheorem{assumption}{\bf Assumption}
\newtheorem{remark}{\bf Remark}
\newtheorem{proposition}{\bf Proposition}
\newtheorem{corollary}{\bf Corollary}

\DeclareMathOperator*{\argmax}{arg\,max}
\def\qed{\hfill $\Box$}

\modulolinenumbers[5]

\begin{document}

\begin{frontmatter}

\title{On the asymptotic of exit problems for controlled Markov diffusion processes with random jumps and vanishing diffusion terms}

\author{Getachew K. Befekadu\corref{mycorrespondingauthor}}
\cortext[mycorrespondingauthor]{Corresponding author}
\ead{gbefekadu@ufl.edu}

\address[rvt1]{NRC/AFRL \& Department of Industrial System Engineering, University of Florida - REEF, 1350 N. Poquito Rd, Shalimar, FL 32579, USA}



\begin{abstract}
In this paper, we study the asymptotic of exit problem for controlled Markov diffusion processes with random jumps and vanishing diffusion terms, where the random jumps are introduced in order to modify the evolution of the controlled diffusions by switching from one mode of dynamics to another. That is, depending on the state-position and state-transition information, the dynamics of the controlled diffusions randomly switches between the different drift and diffusion terms. Here, we specifically investigate the asymptotic exit problem concerning such controlled Markov diffusion processes in two steps: (i) First, for each controlled diffusion model, we look for an admissible Markov control process that minimizes the principal eigenvalue for the corresponding infinitesimal generator with zero Dirichlet boundary conditions -- where such an admissible control process also forces the controlled diffusion process to remain in a given bounded open domain for a longer duration. (ii) Then, using large deviations theory, we determine the exit place and the type of distribution at the exit time for the controlled Markov diffusion processes coupled with random jumps and vanishing diffusion terms. Moreover, the asymptotic results at the exit time also allow us to determine the limiting behavior of the Dirichlet problem for the corresponding system of elliptic partial differential equations containing a small vanishing parameter. 
\end{abstract}

\end{frontmatter}

\linenumbers

\section{Introduction} \label{S1}
In this paper, we are mainly concerned with the asymptotic of exit problem for controlled Markov diffusion processes with random jumps and vanishing diffusion, where the random jumps are introduced in order to modify the evolution of the controlled diffusions by switching from one mode of dynamics to another. That is, depending on the state-position and state-transition information, the dynamics of the controlled diffusions randomly switches between the different drift and diffusion terms. In recent years, some interesting studies on the exit problem for dynamical systems with small random perturbations have been widely reported in literature (see, e.g., \cite{VenFre70}, \cite{FreWe84}, \cite{Kif90} and \cite{DupKu86} in the context of large deviations; see \cite{DupKu89}, \cite{Day86}, \cite{EvaIsh85}, \cite{Zab85}, \cite{Fle78}, \cite{FleTs81} and \cite{BefA15a} in connection with stochastic optimal control problems; and see \cite{Day86} or \cite{MatSc77} via an asymptotic expansions approach). Moreover, we also observe that the exit problem for small random perturbations of dynamical systems with or without random jumps is still a source of interesting problems in engineering, mathematical physics and finance, and elsewhere (e.g., see \cite{BisB11}, \cite{Herr05}, \cite{ImkP06} and the references therein). Note that, in some sense, the rationale behind our framework follows the settings of these papers. In particular, we study the asymptotic exit problem concerning such controlled Markov diffusion processes coupled with random jumps and vanishing diffusion terms in two steps: (i) First, for each controlled diffusion model, we look for an admissible Markov control process that minimizes the principal eigenvalue for the corresponding infinitesimal generator with zero Dirichlet boundary conditions, where such an admissible control process also forces the controlled diffusion process to remain in a given bounded open domain for a longer duration. Moreover, as the diffusion term vanishes, the corresponding deterministic dynamical system is assumed to have a non-empty maximal closed invariant set in the given bounded open domain. (ii) Then, using large deviations theory, we determine the exit place and the type of distribution at the exit time for the controlled Markov diffusion processes coupled with random jumps and vanishing diffusion term. As a consequence of this, the asymptotic results at the exit time also determine the limiting behavior of the Dirichlet problem for the corresponding system of elliptic partial differential equations (PDEs) containing a small vanishing parameter.

The rest of the paper is organized as follows. In Section~\ref{S2}, we present some preliminary results that are useful for our main results. In Section~\ref{S3}, using the basic remarks made in Sections~\ref{S2}, we discuss the action functional for a class of dynamical systems with random jumps and vanishing diffusion terms. In this section, we also briefly discuss the Dirichlet problem for the corresponding system of elliptic PDEs with small vanishing parameter. In Section~\ref{S4}, we provide our main results on the asymptotic estimates for the joint type occupation times and exit distributions for the controlled Markov diffusion processes and the solutions for the corresponding Dirichlet problem with small vanishing parameter. 

\section{Preliminaries} \label{S2}

Consider a family of $d$-dimensional controlled diffusion processes satisfying, in the It\^{o} sense, the following system of stochastic differential equations (SDEs)
\begin{align} 
 dX_{k}^{\epsilon,u_k}(t) = f_k\bigl(X_{k}^{\epsilon, u_k}(t), u_{k}^{\epsilon}(t)\bigr) dt +  \sqrt{\epsilon} \sigma_k\bigl(X_{k}^{\epsilon, u_k} (t)\bigr) d W(t), \notag\\
            X_{k}^{\epsilon,u_k}(0)=x_{0}, \quad  k = 1, 2, \ldots, n, \label{Eq1} 
\end{align}
where
\begin{itemize}
\item $X_{k}^{\epsilon, u_k}(\cdot)$ is an $\mathbb{R}^{d}$-valued controlled diffusion process,
\item $u_{k}^{\epsilon}(\cdot)$ is a $U_k$-valued admissible control process\footnote{Such an admissible control is also assumed to have measurable paths, with $u_{k}^{\epsilon}(t, \omega)\colon [0, \infty) \times \Omega \rightarrow U_k$, for $0 < \epsilon \ll 1$, and satisfies
\begin{align*} 
 \mathbb{E} \left\{\int_0^t {\bigl \vert u_k^{\epsilon}(t)\bigr \vert^2} dt \right\} < \infty, \quad \forall t > 0, \quad  k = 1, 2, \ldots, n. 
\end{align*}} (where $U_k \subset \mathbb{R}^{m_k}$ is compact) which is nonanticipatory and jointly measurable in $(t, \omega) \in [0, \infty) \times \Omega$, on a complete probability space $(\Omega, \mathcal{F}, \mathbb{P})$,
\item $f_k\bigl(\cdot, \cdot\bigr)$, for $k=1,2, \ldots, n$, are bounded Lipschitz continuous functions,
\item $\sigma_k\bigl(\cdot\bigr)$, for $k=1,2, \ldots, n$, are of class $C^2$ functions and bounded together with their derivatives such that $a_k(x) = \sigma_k(x)\,\sigma_k^{T}(x)$, $k =1, 2, \ldots, n$, are uniformly elliptic, i.e.,
\begin{align*}
 a_{min} \vert p \vert^2 \le p \cdot a_k\bigl(x\bigr) p \le a_{max} \vert p \vert^2, \quad x \in \mathbb{R}^{d}, \quad \forall p \in \mathbb{R}^{d}, 
\end{align*}
for some $a_{max} > a_{min} > 0$, 
\item $W(\cdot)$ is a standard Wiener process in $\mathbb{R}^{d}$, and
\item $\epsilon$ is a small parameter, i.e., $0 < \epsilon \ll 1$ that represents the level of random perturbation.
\end{itemize}
Let $\tau_{k}^{\epsilon, u_k}$ be the first exit time corresponding to the diffusion process $X_{k}^{\epsilon,u_k}(t)$ from a given bounded open domain $D \subset \mathbb{R}^{d}$, with a sufficiently smooth boundary $\partial D$, i.e.,
\begin{align}
\tau_{k}^{\epsilon, u_k} = \min \bigl\{ t > 0 \, \bigl\vert \, X_{k}^{\epsilon,u_k}(t) \notin D \bigr\}. \label{Eq2}
\end{align}
As usual, we denote by $\mathbb{P}_{x_0, k}^{\epsilon,u_k} \bigl\{\cdot\bigr\}$ the probability measures in the space of trajectories for the controlled Markov diffusion process $X_{k}^{\epsilon,u_k}(t)$ starting at $x_0 \in D$ and by $\mathbb{E}_{x_0, k}^{\epsilon, u_k} \bigl\{\cdot\bigr\}$ the associated expectation.

In what follows, we consider a particular class of admissible Markov control process $u_k^{\epsilon}$ of the form
\begin{align}
u_k^{\epsilon}(t) = v_k^{\epsilon}\bigl(X_{k}^{\epsilon,v_k}(t)\bigr), \quad \forall t, \quad k =1,2, \ldots, n, \label{Eq3} 
\end{align}
for some measurable function $v_k^{\epsilon} \colon \mathbb{R}^{d} \times \Omega \rightarrow U_k$ and small parameter $\epsilon$ (i.e., $0 < \epsilon \ll 1$). Here, we also identify such an admissible Markov control process, with the map $v_k^{\epsilon}$, from a set of admissible controls $\mathcal{U}_k^{\rm ad}$ for $k=1,2, \ldots, n$.

In reliability theory and other problems arising in engineering, mathematical physics and finance, we often want to confine the controlled diffusion process $X_{k}^{\epsilon,v_k}(t)$ in the given bounded domain $D$ by appropriately choosing admissible controls from the set $\prod_{k=1}^n\mathcal{U}_k^{\rm ad}$ that force the process $X_{k}^{\epsilon,v_k}(t)$ to remain in the bounded domain $D$, {\it at least on the average}, for a longer duration. Note that a standard formulation for such a problem is to maximize the mean exit time of the controlled  diffusion process $X_{k}^{\epsilon,v_k}(t)$ from $D$, i.e.,
\begin{align}
\max_{v_{k}^{\epsilon} \in \mathcal{U}_k^{\rm ad}} \mathbb{E}_{x, k}^{\epsilon, v_k} \left \{\tau_{k}^{\epsilon, v_k} \right\}, \quad x \in D, \label{Eq4} 
\end{align}
or to minimize the exit probability 
\begin{align}
\min_{v_{k}^{\epsilon} \in \mathcal{U}_k^{\rm ad}} \mathbb{P}_{x, k}^{\epsilon, v_k} \left \{\tau_{k}^{\epsilon, v_k} \le T \right\}, \quad x \in D, \label{Eq5} 
\end{align}
for each $k \in \{1, 2, \ldots, n\}$ and small parameter $0 < \epsilon \ll 1$; and we further make a connection to the asymptotic behavior of the exit probabilities as $\epsilon \rightarrow 0$ (i.e., when the diffusion terms vanish).

In general, it is difficult to get effective information about a minimum exit probability and, at the same time, a set of admissible Markov controls in this way. Instead, here we seek to minimize a more natural object, i.e., the exit rate with which the controlled diffusion process $X_{k}^{\epsilon,v_k}(t)$ exits from the given bounded domain $D$; and we further investigate the asymptotic behavior as $\epsilon \rightarrow 0$ and $t \rightarrow \infty$. Note that this leads us to minimize the following family of exit rates 
\begin{align}
\lambda_{k}^{\epsilon, v_k}  = -\limsup_{t \rightarrow \infty} \frac{1}{t} \log \mathbb{P}_{x, k}^{\epsilon, v_k} \left \{\tau_{k}^{\epsilon, v_k} > t \right\}, \,\,\, x \in D, \,\,\, k=1,2, \ldots, n, \label{Eq6} 
\end{align}
w.r.t. some admissible Markov controls $v_{k}^{\epsilon} \in \mathcal{U}_k^{\rm ad}$ for $k =1,2, \ldots, n$.

In order to make clear the above argument, let us define the following family of contraction semigroups
\begin{align}
\mathcal{T}_{t,k}^{\epsilon, v_k} \psi(x) = \mathbb{E}_{x, k}^{\epsilon, v_k} \left \{\psi\bigl(X_{k}^{\epsilon,v_k}(t) \bigr);\, \tau_{k}^{\epsilon, v_k} > t  \right\}, \,\,\, x \in D, \,\,\, k = 1, 2, \ldots, n, \label{Eq7} 
\end{align}
where $\psi \in C_0(\bar{D})$, with $C_0(\bar{D})$ is the space of continuous functions on $\bar{D} \equiv D \cup \partial D$ vanishing on the boundary $\partial D$, and $X_{k}^{\epsilon,v_k}(t)$ satisfies the following SDE
\begin{align} 
 dX_{k}^{\epsilon,v_k}(t) = f_k\bigl(X_{k}^{\epsilon, v_k}(t), v_{k}^{\epsilon}\bigl(X_{k}^{\epsilon, v_k}(t) \bigr)\bigr) dt +  \sqrt{\epsilon} \sigma_k\bigl(X_{k}^{\epsilon, v_k} (t)\bigr) d W(t), \label{Eq8} 
\end{align}
with $x_{k}^{\epsilon,v_k}(0)=x \in D$ and $v_{k}^{\epsilon} \in \mathcal{U}_k^{\rm ad}$ for $k \in \{1,2, \ldots, n\}$. 

Notice that the corresponding infinitesimal generators for the above family of semigroups are extension of the operators, $\mathcal{L}_{k}^{\epsilon, v_k}$ (w.r.t. the admissible Markov controls $v_{k}^{\epsilon} \in \mathcal{U}_k^{\rm ad}$, $k = 1, 2, \ldots, n$), defined on $C^2(\bar{D}) \cap \bigl\{ \psi \, \bigl \vert \, \psi \in C_0(\bar{D}) \,\, \text{and} \,\, \mathcal{L}_{k}^{\epsilon, v_k} \psi \in C_0(\bar{D})\bigr\}$ by
\begin{align} 
\mathcal{L}_{k}^{\epsilon, v_k} \psi(x)= \bigl \langle \bigtriangledown_x \psi(x), f_k(x, v_k^{\epsilon}(x)) \bigr\rangle + \frac{\epsilon}{2} \operatorname{tr}\bigl \{a_k(x)\bigtriangledown_x^2 \psi(x)\bigr\}, \label{Eq9} 
\end{align}
with zero Dirichlet boundary conditions. Moreover, $\mathcal{T}_{t,k}^{\epsilon, v_k}$ is compact for each $k \in \{1,2, \ldots, n\}$ and small (fixed) parameter $0 < \epsilon \ll 1$. Then, the corresponding resolvent operator is compact and, hence, the spectrum $\operatorname{Sp} \bigl\{ \mathcal{L}_{k}^{\epsilon, v_k}\bigr\}$ consists only eigenvalues for $k \in \{1,2, \ldots, n\}$. Furthermore, from Krein-Rutman theorem (e.g., see \cite[Chapter~3, Theorem~3.1]{Pin95}), we also observe that the generator $\mathcal{L}_{k}^{\epsilon, v_k}$, for $k \in \{1,2, \ldots, n\}$, possesses a principal eigenvalue $\lambda_{k}^{\epsilon, v_k}$, which is simple, real and satisfies the following
\begin{align} 
\lambda_{k}^{\epsilon, v_k} = \inf \Bigl \{ \operatorname{Re} \bigl\{ \lambda\bigr\} \, \bigl \vert \, \lambda \in \operatorname{Sp} \bigl\{ \mathcal{L}_{k}^{\epsilon, v_k}\bigr\} \Bigr\}, \,\,\, k =1,2, \ldots, n. \label{Eq10} 
\end{align}
Note that $\lambda_{k}^{\epsilon, v_k}  \in \operatorname{Sp} \bigl\{ \mathcal{L}_{k}^{\epsilon, v_k}\bigr\}$ if and only if $\exp\bigl\{-\lambda_{k}^{\epsilon, v_k}\bigr\} \in \operatorname{Sp} \bigl\{\mathcal{T}_{t,k}^{\epsilon, v_k}\bigr\}$ (e.g., see \cite[Section~3.4]{Paz83}). Hence, since $\mathcal{T}_{t,k}^{\epsilon, v_k} \psi(x) = \mathbb{E}_{x, k}^{\epsilon, v_k} \left \{\psi\bigl(X_{k}^{\epsilon,v_k}(t) \bigr);\, \tau_{k}^{\epsilon, v_k} > t  \right\}$ is compact; and $\Vert \mathcal{T}_{t,k}^{\epsilon, v_k}\Vert < 1$, with $\psi \in C_0(\bar{D})$, $\Vert \psi \Vert < 1$, for each $k \in \{1,2, \ldots, n\}$. Then, it follows that $\lambda_{k}^{\epsilon, v_k} > 0$ and we have the following
\begin{align}
\sup_{\psi \in C_0(\bar{D}),\, \Vert \psi \Vert < 1} \left \Vert \mathcal{T}_{t,k}^{\epsilon, v_k} \psi(x) \right\Vert  = \sup_{x \in D} \mathbb{P}_{x,k}^{\epsilon, v_k} \Bigl \{\tau_{k}^{\epsilon, v_k} > t \Bigr\}, \,\,\,  k =1,2, \ldots, n, \label{Eq11} 
\end{align}
w.r.t. some admissible Markov controls $v_{k}^{\epsilon} \in \mathcal{U}_k^{\rm ad}$ for $k = 1,2, \ldots, n$. Then, a standard argument from \cite{ReSi72} (e.g., see \cite[Chapters~6 and 7]{ReSi72}) also allows us to conclude that
\begin{align}
\lambda_{k}^{\epsilon, v_k} = -\lim_{t \rightarrow \infty} \frac{1}{t} \log \sup_{x \in D} \mathbb{P}_{x, k}^{\epsilon, v_k} \Bigl \{\tau_{k}^{\epsilon, v_k} > t \Bigr\}, \,\,\, x \in D, \,\,\, k \in \{1,2, \ldots, n\} \label{Eq12} 
\end{align}
and which is further equivalent to
\begin{align*}
\lambda_{k}^{\epsilon, v_k} = -\limsup_{t \rightarrow \infty} \frac{1}{t} \log \mathbb{P}_{x, k}^{\epsilon, v_k} \Bigl \{\tau_{k}^{\epsilon, v_k} > t \Bigr\}, \,\,\, x \in D, \,\,\, k \in \{1,2, \ldots, n\}.
\end{align*}
Here, it is worth remarking that the asymptotic behavior for $\mathbb{P}_{x, k}^{\epsilon, v_k} \left \{\tau_{k}^{\epsilon, v_k} > t \right\}$ and $\mathbb{E}_{x, k}^{\epsilon, v_k} \left \{\tau_{k}^{\epsilon, v_k} \right\}$ determine whether the following class of deterministic dynamical systems
\begin{align} 
 \dot x_{k}^{0,u_k}(t) = f_k\bigl(x_{k}^{0, u_k}(t), u_{k}^{0}(t)\bigr), \quad  x_{k}^{0,u_k}(0)=x, \quad  k = 1, 2, \ldots, n, \label{Eq13} 
\end{align}
have a non-empty maximal closed invariant set in $\bar{D}$ or not (e.g., see \cite{Kif81} for additional discussions; and see also Remark~\ref{R1} below). Moreover, under some conditions, the asymptotic behavior for $\mathbb{P}_{x,k}^{\epsilon, v_k} \left \{\tau_{k}^{\epsilon, v_k} > t \right\}$ is also determined by the asymptotics of the principal eigenvalue for $-\mathcal{L}_{k}^{\epsilon, v_k}$, for $k \in \{1,2, \ldots, n\}$ w.r.t. some measurable functions $v_k ^{\epsilon} \colon \mathbb{R}^{d} \times \Omega \rightarrow U_j$, with $u_j^{\epsilon}(t) = v_j^{\epsilon}(x(t))$, $k = 1, 2, \ldots, n$, for {\it almost everywhere} (a.e.) $t \in [0, T]$ as $\epsilon \rightarrow 0$. Moreover, we also observe that the asymptotics of the principal eigenvalue for $-\mathcal{L}_{k}^{\epsilon, v_k}$, for $k \in \{1,2, \ldots, n\}$, are independent of the initial condition $x \in D$. 

Throughout this paper, we assume that there exist a family of measurable functions $\hat{v}_k^{\epsilon}$, with $\hat{u}_k^{\epsilon}(t) = \hat{v}_k^{\epsilon}(x(t))$, for a.e. $t \in [0, T]$, from the set $\mathcal{U}_k^{\rm ad}$, $k = 1, 2, \ldots, n$, rendering a non-empty maximal closed invariant set in $\bar{D}$, as $\epsilon \rightarrow 0$, for the above class of deterministic dynamical systems in \eqref{Eq13}. Later in Section~\ref{S2}, we also require that such a class of deterministic dynamical systems to satisfy additional assumptions (i.e., the large deviations condition of Assumptions~\ref{AS1}; and see also Assumptions~\ref{AS2}).

\begin{remark} \label{R1}
Note that estimating the asymptotic exit rate with which the process $X_{k}^{\epsilon,v_k}(t)$ exits from the given domain $D$ is related to a singularly perturbed eigenvalue problem. For example, the asymptotic behavior for the principal eigenvalue corresponding to the following eigenvalue problem
\begin{align}
\left.\begin{array}{c}
 -\mathcal{L}_k^{\epsilon, v_k}\,\psi_k^{\epsilon,v_k}(x) = \lambda_k^{\epsilon,v_k}\,\psi_k^{\epsilon,v_k}(x)  \quad \text{in} \quad D\\
 \psi_k^{\epsilon,v_k}(x) = 0 \quad \text{on} \quad \partial D, \quad k =1,2, \ldots, n
\end{array}\right\}, \label{Eq14} 
\end{align}
where $\psi_k^{\epsilon,v_k} \in W_{loc}^{2, p} (D) \cap C(\bar{D})$, for $p > 2$, with $\psi_k^{\epsilon,v_k}(x) > 0$ on $D$, has been well studied in the past (e.g., see \cite{Day83} or \cite{DevFr78} in the context of an asymptotic behavior for the principal eigenfunction). Moreover, for a small parameter $0 < \epsilon \ll 1$, if $\hat{v}_k^{\epsilon}(\cdot)$ is the maximizing measurable selector for $\argmax \bigl\{ \mathcal{L}_k^{\epsilon, v_k}\,\psi_k^{\epsilon,v_k}(x, \cdot)\bigr\}$, with $x \in D$ and $k \in \{1,2, \ldots, n\}$. Then, the corresponding principal eigenvalue is given by
\begin{align*}
\lambda_{k}^{\epsilon, \hat{v}_k} = -\limsup_{t \rightarrow \infty} \frac{1}{t} \log \mathbb{P}_{x, k}^{\epsilon, \hat{v}_k} \bigl \{\tau_{k}^{\epsilon, \hat{v}_k} > t \bigr\}, \,\,\, x \in D, \,\,\, k \in \{1,2, \ldots, n\}.
\end{align*}
Here, a few remarks are appropriate:
\begin{enumerate} [(i)]
\item If $\limsup_{\epsilon \rightarrow 0} \lambda_{k}^{\epsilon, \hat{v}_k} < \infty$, for $k = 1, 2, \ldots, n$, and for some $x \in D$, then the maximal closed invariant set for the corresponding deterministic dynamical system $\dot x_{k}^{0,\hat{v}_k}(t) = f_k\bigl(x_{k}^{0,\hat{v}_k}(t), \hat{v}_k^0(x_{k}^{0,\hat{v}_k}(t))\bigr)$, with $x_{k}^{0,\hat{v}_k}(0)=x$, is non-empty in $\bar{D}$ (e.g., see \cite[Theorem~2.1]{Kif81}).
\item If $\limsup_{\epsilon \rightarrow 0} \mathbb{E}_{x, k}^{\epsilon, \hat{v}_k} \bigl \{\tau_{k}^{\epsilon, \hat{v}_k}\bigr\} = \infty$, for $k = 1, 2, \ldots, n$, then the maximal invariant set is non-empty in $\bar{D}$ (cf. \cite[Theorem~3.1]{Kif81}).
\end{enumerate}
\end{remark}

In the present paper, we specifically consider the following controlled Markov process $\bigl(X^{\epsilon,\hat{v}}(t), \zeta^{\epsilon}(t)\bigr)$ in the phase space $\mathbb{R}^d \times \{1,2, \ldots, n\}$
\begin{align} 
 dX^{\epsilon,\hat{v}}(t) = f_{\zeta_t^{\epsilon}} \bigl(X^{\epsilon, \hat{v}}(t), \hat{v}_{{\zeta_t^{\epsilon}}}^{\epsilon}(X^{\epsilon, \hat{v}}(t))\bigr) dt +  \sqrt{\epsilon} \sigma_{{\zeta_t^{\epsilon}}}\bigl(X^{\epsilon,\hat{v}}(t)\bigr) d W(t),  \label{Eq15} 
\end{align}
where
\begin{itemize}
\item $X^{\epsilon,\hat{v}}(t)$ is an $\mathbb{R}^{d}$-valued diffusion process,
\item $\zeta_t^{\epsilon} \equiv \zeta^{\epsilon}(t)$ is a $\{1,2, \ldots, n\}$-valued process such that 
\begin{align*}
 \mathbb{P} \Bigl\{ \zeta^{\epsilon}(t+\triangle) = m \, \bigl\vert \, \zeta^{\epsilon}(t) = k, X_t^{\epsilon,\hat{v}} = x \Bigr\} = \gamma_{km}(x) \triangle + o(\triangle)  \,\,\, \text{as} \,\,\, \triangle \downarrow 0, 
\end{align*}
for $k, m \in \{1,2, \ldots, n\}$ and $k \neq m$, and
\item $\hat{v}_k^{\epsilon}(\cdot) \in \mathcal{U}_k^{\rm ad}$, for $k=1,2, \ldots, n$, with a small parameter $0 < \epsilon \ll 1$, are priorly determined admissible Markov control processes satisfying some additional conditions.\footnote{Here, such admissible Markov control processes are assumed to ensure a maximal closed invariant set in $\bar{D}$ for the corresponding class of deterministic dynamical systems in \eqref{Eq13} as the small parameter goes to zero, i.e., when $\epsilon \rightarrow 0$.}
\end{itemize}

\begin{remark} \label{R2}
Note that, in equation~\eqref{Eq15}, the random jumps, depending on the state-position and state-transition information, modify the dynamics of the controlled diffusions by switching randomly from one mode of operation to another.
\end{remark}

Here, we also assume that the transition coefficients $\gamma_{km}(x)$, for $x \in \mathbb{R}^{d}$, are positive and Lipschitz continuous. Moreover, under these conditions (e.g., see \cite{EizF90} and \cite[Chapter~2]{Fre96}), there exists a unique vector $\bar{\omega}(x) = \bigl(\omega_1(x), \omega_2(x), \ldots, \omega_n(x)\bigr)$ such that
\begin{align}
 \omega_k(x) > 0, \quad \sum\nolimits_{k = 1}^{n} \omega_k(x) = 1 \quad \text{and} \quad \bar{\omega}(x) \Gamma(x) = 0, \label{Eq16} 
\end{align}
where $\Gamma(x) = \bigl(\Gamma_{km}(x)\bigr)$ is an $n \times n$ matrix and
\begin{align}
\Gamma_{km}(x) = \left \{ \begin{array}{l}
  \gamma_{km}(x) \quad\qquad \quad\quad \quad \text{for} \quad k \neq m\\
 - \sum\nolimits_{j: j \neq k} \gamma_{kj}(x) \quad \quad\, \text{for}  \quad k = m,
\end{array} \right. \label{Eq17} 
\end{align}
with $k, m \in \{1, 2, \ldots, n\}$.

Denote by $\mathbb{P}_{x_0, k}^{\epsilon,\hat{v}} \bigl\{\cdot\bigr\}$ the probability measures in the space of trajectories of the process $(X^{\epsilon,\hat{v}}(t), \zeta^{\epsilon}(t))$ and by $\mathbb{E}_{x_0, k}^{\epsilon,\hat{v}} \bigl\{\cdot\bigr\}$ the associated expectation. Define the occupation time $r^{\epsilon,\hat{v}}(t)$ for the component $\zeta^{\epsilon}(t)$ as 
\begin{align}
 r^{\epsilon,\hat{v}}(t) &= \bigl(r_{1}^{\epsilon,\hat{v}}(t), r_{2}^{\epsilon,\hat{v}}(t), \ldots, r_{n}^{\epsilon,\hat{v}}(t)\bigr) \notag \\
                      &= \biggl(\int_0^t \chi_1\bigl(\zeta^{\epsilon}(s)\bigr) ds, \int_0^t \chi_2\bigl(\zeta^{\epsilon}(s)\bigr) ds, \ldots, \int_0^t \chi_n\bigl(\zeta^{\epsilon}(s)\bigr) ds \biggr),
                        \label{Eq18} 
\end{align}
where $r_{k}^{\epsilon,\hat{v}}(t) = \int_0^t \chi_k\bigl(\zeta^{\epsilon}(s)\bigr) ds \in \mathbb{R}_{+}$, for any $t \in [0, T]$, and $\chi_k$ is the indicator function of the singleton set $\{k\}$, with $k=1,2, \ldots, n$. Then, we specifically study the controlled Markov process  $(X^{\epsilon,\hat{v}}(t), \zeta^{\epsilon}(t))$ and the occupation time $r^{\epsilon,\hat{v}}(t)$; and we further investigate the component process $X^{\epsilon,\hat{v}}(t)$ as a result of small random perturbations of the following average dynamical system
\begin{align}
 \dot{x}(t) &= \sum\nolimits_{k=1}^n \omega_k(x(t)) f_k \bigl(x(t), \hat{v}_k^{0}(x(t))\bigr)  \notag \\
                &\triangleq f_{\rm av}^{\hat{v}} \bigl(x(t)\bigr), \quad \quad x(0) = x_0 \in \mathbb{R}^d,  \label{Eq19}  
\end{align}
where $\hat{v}_k^{0}$, $k=1,2, \ldots, n$, are some non-anticipatory maps (with $\hat{u}_k^{0}(t) = \hat{v}_k^{0}(x(t))$, for {\rm a.e.} $t \in [0, T]$) that allows us to prove large deviations results for the joint type occupation times and positions as $\epsilon \rightarrow 0$ and study the exit probabilities for such a family of processes.

On the other hand, the infinitesimal generator $\mathcal{L}^{\epsilon,\hat{v}}$ of the process $(X^{\epsilon,\hat{v}}(t), \zeta^{\epsilon}(t))$ acting on smooth functions (smooth in $x \in \mathbb{R}^d$) is given by
\begin{align}
 \mathcal{L}^{\epsilon,\hat{v}}\psi_k(x)= \mathcal{L}_k^{\epsilon,\hat{v}_k}\psi_k(x) + \sum\nolimits_{j=1}^n \gamma_{kj}(x) \bigl[\psi_j(x) - \psi_k(x)\bigr], \label{Eq20}
\end{align}
where
\begin{align}
\mathcal{L}_{k}^{\epsilon,\hat{v}_k} \psi_k(x) = \bigl \langle \bigtriangledown_x \psi_k(x), f_k(x, \hat{v}_k(x)) \bigr\rangle + \frac{\epsilon}{2} \operatorname{tr}\bigl \{a_k(x)\bigtriangledown_x^2 \psi_k(x) \bigr\}. \label{Eq21}
\end{align}

\begin{remark} \label{R3}
Later in Section~\ref{S4}, we also establish a connection between those controlled Markov diffusion processes corresponding to the infinitesimal generators $\mathcal{L}_k^{\epsilon,\hat{v}}$, for $k=1,2, \ldots, n$, and that of the occupation time $r^{\epsilon,\hat{v}}(t)$, with $t \in [0, T]$, for the component process $\zeta^{\epsilon}(t)$ (e.g., see Proposition~\ref{P3}).
\end{remark}

Finally, note that the controlled Markov process $(X^{\epsilon,\hat{v}}(t), \zeta^{\epsilon}(t))$ is closely connected with the following Dirichlet problem corresponding to a system of elliptic PDEs that satisfies the maximum principle (e.g., see \cite[Chapter~3, Section~8]{ProW84} for the application of maximum principle for classical Dirichlet problems),
\begin{align}
\left \{ \begin{array}{l}
  \mathcal{L}_k^{\epsilon,\hat{v}} \psi_k^{\epsilon,\hat{v}}(x) + \sum\nolimits_{j=1}^n \gamma_{kj}(x) \bigl[\psi_j^{\epsilon,\hat{v}}(x) - \psi_k^{\epsilon,\hat{v}}(x)\bigr] = 0, \quad x \in D,\\
 \psi_k^{\epsilon,\hat{v}}(x)\vert_{\partial D} = g_k(x), \quad k=1,2, \ldots, n,
\end{array} \right. \label{Eq22}
\end{align}
where we can also study the limiting behavior for the solution of the above Dirichlet problem as the small parameter vanishes, i.e., when $\epsilon \rightarrow 0$. Here, we remark that the interplay between the small diffusion and the jumps $\zeta$-component leads to the situation -- where $g_k(x)$, for $k=1, 2, \ldots, n$, will influence the $\lim_{\epsilon \downarrow 0} \psi_k^{\epsilon,\hat{v}}(x)$. 

\begin{remark} \label{R4}
Note that, when $\epsilon=0$, the controlled Markov diffusion process $(X^{\epsilon,\hat{v}}(t), \zeta^{\epsilon}(t), \mathbb{P}_{x_0, k}^{\epsilon,\hat{v}})$ turns into a random process $(X^{0,\hat{v}}(t), \zeta^{0}(t), \mathbb{P}_{x_0, k}^{0,\hat{v}})$, corresponding to the following system of first-order equations
\begin{align*}
 \mathcal{L}^{0,\hat{v}}\psi_k^{0,\hat{v}}(x) &= \bigl \langle \bigtriangledown_x \psi_k^{0,\hat{v}}(x), f_k(x, \hat{v}_k(x)) \bigr\rangle + \sum\nolimits_{j=1}^n \gamma_{kj}(x) \bigl[\psi_j^{0,\hat{v}}(x) - \psi_k^{0,\hat{v}}(x)\bigr]\\
& = 0, \quad  x \in D, \quad k=1,2, \ldots, n, 
\end{align*}
where such a random process $(X^{0,\hat{v}}(t), \zeta^{0}(t), \mathbb{P}_{x_0, k}^{0,\hat{v}})$ consists of deterministic motion along the fields $f_k(x, \hat{v}_k^{0}(x))$, $k=1,2, \ldots, n$, and the random switching from one field to another (see \cite[Chapter~2]{Fre96}). Moreover, the random exit time $\tau_{D}^{0,\hat{v}} = \bigl\{t > 0 \,\bigl\vert \, X^{0,\hat{v}}(t) \notin \partial D \bigr\}$, starting from $x_0 \in D$, with a random variable $k \in \{1,2, \ldots, n\}$, further satisfies the following
\begin{align*}
\mathbb{P}_{x_0, k}^{0,\hat{v}} \bigl\{ \tau_{D}^{0,\hat{v}} > t \bigr\} \le \exp\bigl\{-t \lambda^{0,\hat{v}_k} \bigr\}, 
\end{align*}
where $\lambda^{0,\hat{v}_k}$ is a principal eigenvalue for the following problem
\begin{align*}
\left.\begin{array}{c}
 -\mathcal{L}^{0,\hat{v}}\psi_k^{0,\hat{v}_k}(x) = \lambda^{0,\hat{v}_k}\,\psi_k^{0,\hat{v}_k}(x)  \quad \text{in} \quad D\\
 \psi_k^{0,\hat{v}_k}(x) = 0 \quad \text{on} \quad \partial D, \quad k =1,2, \ldots, n
\end{array}\right\},
\end{align*}
with $\psi_k^{0,\hat{v}_k} \in W_{loc}^{2, p} (D) \cap C(\bar{D})$, for $p > 2$, with $\psi_k^{0,\hat{v}_k}(x) > 0$ on $D$. Note such information about $\lambda^{0,\hat{v}_k}$, $k =1,2, \ldots, n$, also agree with our argument as remarked at the beginning of this section (cf. Remark~\ref{R1}).
\end{remark}
In the following sections, we further study the limiting behavior for the solution of the corresponding Dirichlet problem in \eqref{Eq22} in two steps: (i) the first step is related with the asymptotic exit problem for the component $X^{\epsilon,\hat{v}}(t)$ from the domain $D$, where such an exit problem can be addressed by determining the action functional for the family of processes $X^{\epsilon,\hat{v}}(t)$ as $\epsilon \rightarrow 0$, and (ii) the second step is related with determining the position of the component $\zeta^{\epsilon}(t)$ at the random time $\tau_{D}^{\epsilon,\hat{v}} = \bigl\{t > 0 \,\bigl\vert \, X^{\epsilon,\hat{v}}(t) \notin D \bigr\}$.

\section{Action functional for the family $(X^{\epsilon,\hat{v}}(t), r^{\epsilon,\hat{v}}(t))$ as $\epsilon \rightarrow 0$} \label{S3}
In this section, we provide some preliminary results that are concerned with the action functional for the family of processes $(X^{\epsilon,\hat{v}}(t), r^{\epsilon,\hat{v}}(t))$ as $\epsilon$ tends to zero (e.g., see \cite{VenFre70} and \cite[Chapter~3]{FreWe84} for additional discussions on the notion of action functional for dynamical systems with random perturbations). Before stating these results, we need some notations. Let $\lambda(x, p, \alpha)$ be the principal eigenvalue of the matrix $\bigl(H_{km}(x, p,\alpha)\bigr)$, $x, p \in \mathbb{R}^d$, $k, m \in \{1,2, \ldots, n\}$, $\alpha = (\alpha_1, \alpha_2, \ldots, \alpha_n) \in \mathbb{R}^n$
\begin{align}
H_{km}(x, p,\alpha) = \left \{ \begin{array}{l}
  \Gamma_{km}, \hspace{2.07in} \text{if} \quad m \neq k\\
  \bigl[p \cdot a_k(x) p/2 + p \cdot \tilde{f}_k^{\hat{v}}(x) + \alpha_k \bigr] + \Gamma_{kk}, \quad \text{if} \quad m = k,
\end{array} \right. \label{Eq23}
\end{align}
where $\tilde{f}_k^{\hat{v}}(x) = f_k(x, \hat{v}_k^{0}(x))$ for $k =1, 2, \ldots, n$.\footnote{Here, the admissible maps $\hat{v}_k^{0}(x)$, $k=1, 2, \ldots, n$, are assumed to be priorly determined (cf. Section~\ref{S2}).}

Note that $\lambda(x, p, \alpha)$ is convex in $(p, \alpha)$ and its Legendre transform  in $(p, \alpha)$ is given by
\begin{align}
\eta(x, q, \beta) =\sup_{p \in \mathbb{R}^d,\, \alpha \in \mathbb{R}^n}  \bigl[ q \cdot p + \beta \cdot \alpha - \lambda(x, p, \alpha) \bigr], \quad x, q \in \mathbb{R}^d, \,\,\beta \in \mathbb{R}^n. \label{Eq24}
\end{align}

Let $C(\mathbb{R}^d)$ be the space of continuous functions: $[0, T] \rightarrow \mathbb{R}^d$ and 
\begin{align}
C_{+}(\mathbb{R}^n) = \Bigl\{\mu=(\mu_1, \mu_2, \ldots, \mu_n) \,\bigl \vert \,\mu \in C(\mathbb{R}^n), \quad  \mu_k(0)=0, \,\,  \,\, 1 \le k \le n, \quad  & \notag \\
  \mu_k(t)  \,\, \text{is non-decreasing and} \,\, \sum\nolimits_{k=1}^n \mu_k(t) = t, \,\, t \in [0, T] \Bigr\}.& \label{Eq25}
\end{align}

Let $T > 0$ be fixed and define the following functional
\begin{align}
S_{0T}(\varphi, \mu) = \left \{ \begin{array}{l}
  \int_0^T \eta\bigl(\varphi(s), \dot{\varphi}(s), \dot{\mu}(s)\bigr) ds, \quad \text{if} \,\, \varphi \in C(\mathbb{R}^d) \,\, \text{and} \,\, \mu \in C_{+}(\mathbb{R}^n) \\
  \hspace{1.7 in} \text{are absolutely continuous (a.c.)}\\
 +\infty \hspace{1.4 in}\text{otherwise}.
\end{array} \right. \label{Eq26}
\end{align}

Suppose that the diffusion and transition coefficients, i.e., $\sigma_k$ and $\gamma_{km}$, for all $k, m \in \{1,2, \ldots, n\}$, satisfy the Lipschitz continuous and positive-Lipschitz continuous conditions, respectively. Then, we have the following result (w.r.t. the admissible Markov controls $\hat{v}_{k}^{\epsilon} \in \mathcal{U}_k^{\rm ad}$, $k=1, 2, \ldots, n$).
\begin{proposition} \label{P1}
The functional $\epsilon^{-1} S_{0T}$ is the action functional for the family of processes $(X^{\epsilon,\hat{v}}(t), r^{\epsilon,\hat{v}}(t))$ as $\epsilon \rightarrow 0$ in the uniform topology. The rate function, i.e.,  functional $S_{0T}$, is nonnegative and equal to zero only when $\dot{\varphi}(t) = f_{\rm av}^{\hat{v}} \bigl(\varphi(t)\bigr)$ and $\dot{\mu}(t) = \bar{\omega}(\varphi(t))$, for $t \in [0, T]$.
\end{proposition}

In Section~\ref{S3}, we will give a proof for the above proposition after obtaining additional supporting results.

Further, let us denote by $\Psi_t^{x_0}$ the integral curve of the vector field $\dot{x}(t)=f_{\rm av}^{\hat{v}} \bigl(x(t)\bigr)$, with $\hat{v}_k^{0}(x)$, $k=1, 2, \ldots, n$, starting from the point $x(0)=x_0$ (i.e., $\dot{\Psi}_t^{x_0} = f_{\rm av}^{\hat{v}} \bigl(\Psi_t^{x_0}\bigr)$, with $\Psi_0^{x_0} = x_0$). Then, define
\begin{align}
\rho(x, q) =\sup_{p \in \mathbb{R}^d}  \bigl[ q \cdot p - \lambda(x, p, 0) \bigr], \quad x, q \in \mathbb{R}^d  \label{Eq27} 
\end{align}
and
\begin{align}
I_{0T}(\varphi) = \left \{ \begin{array}{l}
  \int_0^T \rho\bigl(\varphi(s), \dot{\varphi}(s)\bigr) ds, \quad \quad \text{if} \,\, \varphi \in C(\mathbb{R}^d) \,\, \text{is a.c.}\\
 +\infty \hspace{1.2 in}\text{otherwise}.
\end{array} \right.  \label{Eq28} 
\end{align}
Taking into account the involution property of the Legendre transform, then we have the following
\begin{align}
\eta(x, q, \beta) &= -\sup_{\beta \in \mathbb{R}^n}  \bigl[ - 0 \cdot \beta - \eta(x, q, \beta)  \bigr]  \notag\\
                         &= \sup_{p \in \mathbb{R}^d}  \bigl[ q \cdot p - \lambda(x, p, 0) \bigr] \notag \\
                         &= \rho(x, q).  \label{Eq29} 
\end{align}

Next, we have the following result which is a direct consequence of the contraction principle (see also \cite[Chapter~5, pp.~117--124]{FreWe84}).
\begin{corollary}  \label{C1}
The functional $\epsilon^{-1} I_{0T}$ is the action functional for the family processes $X^{\epsilon,\hat{v}}(t)$ as $\epsilon \rightarrow 0$ in the uniform topology. Then, the rate function $I_{0T}$ is equal to zero only when $\varphi_t = \Psi_t^{x_0}$ and $\Psi_0^{x_0}=x_0$.
\end{corollary}
Let $\bar{n}(y)$ be a unit vector normal to the boundary $\partial D$ at $y \in \partial D$. Furthermore, we assume that the average dynamical system $f_{\rm av}^{\hat{v}}(x)$ satisfies following large deviations condition.

\begin{assumption} [Large deviations condition] \label{AS1}
The vector field $f_{\rm av}^{\hat{v}}(y)$ points outward from the boundary $\partial D$, i.e., $\bigl\langle f_{\rm av}^{\hat{v}}(y), \bar{n}(y)\bigr\rangle < 0$ for any $y \in \partial D$. Moreover, if the vector field $f_{\rm av}^{\hat{v}}(x)$ has a unique stationary point at $x_0 \in D$. Then, the function
\begin{align}
V(y) = \inf \Bigl\{ I_{0T}(\varphi) \, \bigl\vert \, \varphi(0) = x_0, \,\, T > 0, \,\, \varphi(T) = y \,\,\, \text{for} \,\,\, y \in \partial D \Bigr\}  \label{Eq30} 
\end{align}
attains its unique minimum at $\bar{y}_0 \in \partial D$, i.e., $V(\bar{y}_0) < V(y)$ for any $y \in \partial D$.
\end{assumption}

\begin{assumption} \label{AS2}
There exists $k_0$, with $k_0 \in \{1,2, \ldots, n\}$, such that at the point $\bar{y}_0 \in \partial D$, defined above in Assumption~\ref{AS1}, then the following generic inequalities hold
\begin{align}
 \bigl\langle \tilde{f}_{k_0}^{\hat{v}}(\bar{y}_0), \bar{n}(\bar{y}_0)\bigr\rangle > \bigl\langle \tilde{f}_k^{\hat{v}} (\bar{y}_0), \bar{n}(\bar{y}_0)\bigr\rangle, \,\, k_0 \in \{1,2, \ldots, n\}, \,\, k \neq k_0,  \label{Eq31} 
\end{align}
with $\tilde{f}_k^{\hat{v}}(x) = f_k(x, \hat{v}_k^{0}(x))$, for the non-anticipatory measurable maps $\hat{v}_k^{0}$, $k \in \{1, 2, \ldots, n\}$. 
\end{assumption}

Let $\tau_{D}^{\epsilon,\hat{v}}$ be the first exit time for the component $X^{\epsilon,\hat{v}}(t)$ from $D \times \{1,2, \ldots, n\}$, i.e.,
\begin{align}
\tau_{D}^{\epsilon,\hat{v}} = \min \Bigl\{ t > 0 \, \bigl\vert \, X^{\epsilon,\hat{v}}(t) \notin D \Bigr\}. \label{Eq32}
\end{align}
Then, we can study the limiting distribution of $(X^{\epsilon,\hat{v}}(\tau_{D}^{\epsilon,\hat{v}}), \zeta^{\epsilon}(\tau_{D}^{\epsilon,\hat{v}}))$ as $\epsilon \rightarrow 0$. Moreover, this distribution also determines the limiting behavior for the solution of the Dirichlet problem in \eqref{Eq22}, when the small parameter vanishes, i.e., as $\epsilon \rightarrow 0$.

\section{Main results} \label{S4}
In this section, we present our main results that establish a connection between the asymptotic exit probability problem for $(X^{\epsilon,\hat{v}}(t), \zeta^{\epsilon}(t))$ from $D \times \{1, 2, \ldots, n\}$ and that of the limiting behavior for the solutions of the Dirichlet problem in \eqref{Eq22}, as $\epsilon \rightarrow 0$. Note that if Assumption~\ref{AS1} holds true (i.e., the large deviations condition), then the exit problem for the component $X^{\epsilon,\hat{v}}(t)$ from the given bounded open domain $D$ is equivalent to determining the action functional for the family of processes $X^{\epsilon,\hat{v}}(t)$ as $\epsilon \rightarrow 0$ and that of the exact exit position for the component process $\zeta^{\epsilon}(t)$ at the random time $\tau_{D}^{\epsilon,\hat{v}} = \bigl\{t > 0 \,\bigl\vert \, X^{\epsilon,\hat{v}}(t) \notin \partial D \bigr\}$.

Then, we have our first result concerning the asymptotic estimates for the joint type occupation times and the exit positions.
 \begin{proposition} \label{P2}
Let the diffusion matrices $a_k(x)$ and the transition coefficients $\gamma_{km}(x)$ be Lipschitz continuous and let $a_k(x)$ be uniformly elliptic and $\gamma_{km}(x) > 0$ for $x \in \bar{D}$, $k, m \in \{1,2, \ldots, n\}$, with $k \neq m$. If Assumption~\ref{AS1} holds true, i.e., the large deviations condition. Then, we have
\begin{align}
  \lim_{\epsilon \rightarrow 0} \mathbb{P}_{x_0, k}^{\epsilon,\hat{v}} \bigl \{ \vert X^{\epsilon,\hat{v}}(\tau_{D}^{\epsilon,\hat{v}}) - \bar{y}_0 \vert > \delta \bigr\} = 0, \label{Eq33}
\end{align}
for any $\delta > 0$, $1 \le k \le n$, uniformly in $x_0 \in \hat{D}$ for any compact $\hat{D} \subset D$. Furthermore, if Assumption~\ref{AS2} is satisfied, then we have the following
\begin{align}
  \lim_{\epsilon \rightarrow 0} \mathbb{P}_{x_0, k}^{\epsilon, \hat{v}} \bigl \{\zeta^{\epsilon}(\tau_{D}^{\epsilon,\hat{v}}) = k_0 \, \vert \, \tau_{D}^{\epsilon,\hat{v}} < \infty \bigr\} = 1, \label{Eq34}
\end{align}
for  $1 \le k \le n$ and $x_0 \in \hat{D} \subset D$.
\end{proposition}

Let us establish the following results (i.e., Propositions~\ref{P3} and \ref{P4}) that are useful for proving Proposition~\ref{P2}.

\begin{proposition} \label{P3}
Suppose that the functions $\tilde{f}_k^{\hat{v}}(x)$, $\sigma_k(x)$ and $\gamma_{km}(x)$, for $k, m \in \{1, 2, \ldots, n\}$, are independent of the position variable $x$ (i.e., $f_k^{\hat{v}}(x)$, $\sigma_k(x)$ and $\gamma_{km}(x)$, for $k, m \in \{1, 2, \ldots, n\}$, are constants). Then, the statement in Proposition~\ref{P1} holds true. 
\end{proposition}

In what follows, we consider a variation $\mathbb{Q}_{x_0, k}^{\epsilon,\hat{v}}$ of $\mathbb{P}_{x_0, k}^{\epsilon,\hat{v}}$ that is governed by the same initial value and the evolution except that $\gamma_{km}(t)$ depends on time $t$ rather than state position $x$. Then, the corresponding Legendre transform for $\mathbb{Q}_{x_0, k}^{\epsilon,\hat{v}}$ in $(p, \alpha)$ is given by
\begin{align}
\hat{\eta}(t, x, q, \beta) =\sup_{p \in \mathbb{R}^d,\, \alpha \in \mathbb{R}^n}  \bigl[ q \cdot p + \beta \cdot \alpha - \hat{\lambda}(t, x, p, \alpha) \bigr], \,\,\, x, q \in \mathbb{R}^d, \,\,\, \beta \in \mathbb{R}^n,  \label{Eq43} 
\end{align}
where the principal eigenvalue $\hat{\lambda}(t, x, p, \alpha)$ is associated with the following matrix $\bigl(\hat{H}_{km}(t, x, p,\alpha)\bigr)$, for $k, m \in \{1,2, \ldots, n\}$, with\footnote{Notice that $\Gamma(t) = \bigl(\Gamma_{km}(t)\bigr)$ is an $n \times n$ matrix (cf. Equations~\eqref{Eq17} and \eqref{Eq23}), with 
\begin{align*}
 \Gamma_{km} (t)= \left \{ \begin{array}{l}
  \gamma_{km}(t), \quad \quad\quad \quad \quad \text{if} \quad m \neq k\\
  - \sum\nolimits_{j: j \neq k} \gamma_{kj}(t), \quad \text{if} \quad m = k
\end{array} \right. 
\end{align*}
and $k, m \in \{1,2, \ldots, n\}.$}

\begin{align}
\hat{H}_{km}(t, x,p,\alpha) = \left \{ \begin{array}{l}
  \gamma_{km}(t), \hspace{2.52in} \text{if} \quad m \neq k\\
  \bigl[p \cdot a_k(x) p/2 + p \cdot \tilde{f}_k^{\hat{v}}(x) + \alpha_k \bigr]  - \sum\nolimits_{j: j \neq k} \gamma_{kj}(t), \quad \text{if} \quad m = k.
\end{array} \right.  \label{Eq44} 
\end{align}

Let $T > 0$ be fixed and define
\begin{align}
\hat{S}_{0T}(\varphi, \mu) = \left \{ \begin{array}{l}
  \int_0^T \hat{\eta}\bigl(s, \varphi(s), \dot{\varphi}(s), \dot{\mu}(s)\bigr) ds, \,\,\, \text{if} \,\, \varphi \in C(\mathbb{R}^d) \,\, \text{and} \,\, \mu \in C_{+}(\mathbb{R}^n)\,\, \text{a.c.}\\
 +\infty \hspace{1.45 in}\text{otherwise}.
\end{array} \right.  \label{Eq45} 
\end{align}

Then, we have the following result.

\begin{proposition} \label{P4}
The action functional for the family of processes $(X^{\epsilon,\hat{v}}(t), r^{\epsilon,\hat{v}}(t))$ w.r.t. $\mathbb{Q}_{x_0, k}^{\epsilon,\hat{v}}$ as $\epsilon \rightarrow 0$ is ${\epsilon}^{-1}\hat{S}_{0T}(\varphi, \mu)$ in the uniform topology.
\end{proposition}

{\em Proof of Proposition~\ref{P2}}: 
Let $\hat{D}_{\delta}$ and $\hat{D}_{2\delta}$ be $\delta$ and $2\delta$-neighborhoods of the compact set $\hat{D} \subset D$, with sufficiently smooth boundaries $\partial \hat{D}_{\delta}$ and $\partial \hat{D}_{2\delta}$, respectively. Then, the state-trajectories $X^{\epsilon,v}(t)$, starting from any $x \in D$, $k \in \{1,2, \ldots, n\}$, hit $\partial \hat{D}_{\delta}$ before $\partial D$ with probability close to one as $\epsilon$ is small enough. This follows from Assumption~\ref{AS1}. Hence, taking into account the strong Markov property of the process $(X^{\epsilon,\hat{v}}(t), \zeta^{\epsilon}(t), \mathbb{P}_{x,k}^{\epsilon,\hat{v}})$, it is sufficient to prove Proposition~\ref{P2} for $x \in \partial \hat{D}_{\delta}$, $k \in \{1,2, \ldots, n\}$.

Define the following Markov times $\theta_0 < \tau_1 < \theta_1 < \cdots < \tau_{\ell} < \theta_{\ell} \cdots$ as follows 
\begin{align*}
\left. \begin{array}{l}
 \theta_0 = \min \bigl\{t > 0 \,\vert \, X^{\epsilon,\hat{v}}(t) \in\partial \hat{D}_{2\delta} \bigr\}\\
 \tau_1 = \min \bigl\{t > \theta_0 \,\vert \, X^{\epsilon,\hat{v}}(t) \in\partial \hat{D}_{\delta} \cup \partial D \bigr\}\\
 \theta_1 = \min \bigl\{t > \tau_1 \,\vert \, X^{\epsilon,\hat{v}}(t) \in\partial \hat{D}_{2\delta} \bigr\}\\
  \hspace{0.75 in} \cdots \\
 \tau_{\ell+1} = \min \bigl\{t > \theta_{\ell} \,\vert \, X^{\epsilon,\hat{v}}(t) \in\partial \hat{D}_{\delta} \cup \partial D \bigr\}\\
 \theta_{\ell+1} = \min \bigl\{t > \tau_{\ell} \,\vert \, X^{\epsilon,\hat{v}}(t) \in\partial \hat{D}_{2\delta} \bigr\}\\
  \hspace{0.75 in } \cdots 
\end{array} \right.
\end{align*}
Next, let us define a Markov chain $(\hat{X}_{\ell}^{\epsilon,\hat{v}}, \hat{\zeta}_{\ell}^{\epsilon})$ in the phase space $\bigl\{\hat{D}_{\delta} \cup \partial D\bigr\} \times \bigl\{1,2, \ldots, n\bigr\}$ with $\hat{X}_{\ell}^{\epsilon,\hat{v}} = X^{\epsilon,\hat{v}}(\tau_{\ell})$ and $\hat{\zeta}_{\ell}^{\epsilon} = \zeta^{\epsilon}(\tau_{\ell})$. Note that the first exit of $X^{\epsilon,\hat{v}}(t)$ from the domain $D$ occurs, when the component $\hat{X}_{\ell}^{\epsilon,\hat{v}}$ of the chain first time belongs to $\partial D$. Then, using the large deviations estimate for the family of processes $(X^{\epsilon,\hat{v}}(t), \zeta^{\epsilon}(t), \mathbb{P}_{x,k}^{\epsilon,\hat{v}})$, as $\epsilon \rightarrow 0$, we can show, in the standard way (e.g., see \cite[Chapter~4]{FreWe84}), that $\hat{X}_{\ell}^{\epsilon,\hat{v}}$ starting from any $x \in \partial \hat{D}_{\delta}$ and $k \in \{1,2, \ldots, n\}$ reaches $\partial D$ for the first time to a small neighborhood of the point $\bar{y}_0 \in \partial D$, introduced in Assumption~\ref{AS2}, with probability close to one as both parameters $\epsilon$ and $\delta$ are small enough, which implies the first statement of Proposition~\ref{P2}.

In order to prove the second statement, we use the fact that the extremal of the variational problem
\begin{align}
\inf \Bigl\{ I_{0T}(\varphi) \, \bigl\vert \, \varphi(0) \in \hat{D}, \,\,\, \varphi(T) \in \partial D, \,\, T > 0 \Bigr\}  \label{Eq56} 
\end{align}
spends in $\delta$-neighborhood $\mathscr{E}_{\delta} = \bigl\{x \in D \, \vert \, \rho(x, \partial D) < \delta \bigr\}$ of $\partial D$ a time of order $\delta$ as $\delta \rightarrow 0$. Note that, with probability close to one as $\delta$ is small, the second component $\zeta^{\epsilon}(t)$ has no jumps during this time; and, hence, $X^{\epsilon,v}(t)$ hits the boundary for the value of the second coordinate $\zeta^{\epsilon}(t)$ such that the transition of  $X^{\epsilon,v}(t)$ from $\mathscr{E}_{\delta} \setminus \partial D$ to $\partial D$ is easiest transition, when the second component is equal to $k_0$ defined in Assumption~\ref{AS2}. This completes the proof.\qed

Note that the limiting distributions of $(X^{\epsilon,\hat{v}}(\tau_{D}^{\epsilon,\hat{v}}), \zeta^{\epsilon}(\tau_{D}^{\epsilon,\hat{v}}))$ as $\epsilon \rightarrow 0$ also determine the limiting behavior for the solutions of the corresponding Dirichlet problem in \eqref{Eq22}, where such a connection is further established using the following result.
\begin{proposition} \label{P5}
Suppose Assumptions~\ref{AS1} and \ref{AS2} hold true. Then, we have
\begin{align}
\lim_{\epsilon \rightarrow 0} \psi_k^{\epsilon,\hat{v}}(x) = g_{k_0} (\bar{y}_0), \,\, 1 \le k \le n, \label{Eq57}
\end{align}
uniformly in $x \in \hat{D} \subset D$, where $k_0 \in \{1, 2, \ldots, n\}$, $\bar{y}_0 \in \partial D$ (cf. Assumptions~\ref{AS2}) and  $\psi_k^{\epsilon,\hat{v}}(x)$ is the solution for the Dirichlet problem in \eqref{Eq22}.
\end{proposition}

{\em Proof}: 
The proof easily follows from Proposition~\ref{P2} and the stochastic representation of
\begin{align}
\psi_k^{\epsilon,\hat{v}}(x) = \mathbb{E}_{x,k}^{\epsilon,\hat{v}} \Bigl\{ \psi_{\zeta^{\epsilon}(\tau_D^{\epsilon,\hat{v}})}^{\epsilon,\hat{v}}(X^{\epsilon,\hat{v}}(\tau_{D}^{\epsilon,\hat{v}})) \Bigr\},  \label{Eq58} 
\end{align}
which is uniformly in $x \in \hat{D} \subset D$ (cf. \cite[Theorem~3]{EizF90}). Furthermore, notice that
\begin{align}
\mathbb{P}_{x,k}^{\epsilon,\hat{v}} \Bigl\{ \tau_D^{\epsilon,\hat{v}} < \infty \Bigr\} = 1, \,\,\, \text{for any} \,\, x \in \hat{D}, \,\ k \in \{1,2, \ldots, n\}.  \label{Eq59} 
\end{align}
Thus, taking into account $\lim_{\epsilon \rightarrow 0} \mathbb{P}_{x,k}^{\epsilon,\hat{v}} \bigl\{ \zeta^{\epsilon}(\tau_D^{\epsilon,\hat{v}}) = k_0 \bigr\} = 1$ (see Assumption~\ref{AS2}), we have
\begin{align}
\lim_{\epsilon \rightarrow 0} \psi_k^{\epsilon,\hat{v}}(x) &= \lim_{\epsilon \rightarrow 0} \mathbb{E}_{x,k}^{\epsilon,\hat{v}} \Bigl\{ \psi_{\zeta^{\epsilon}(\tau_D^{\epsilon,\hat{v}})}^{\epsilon,\hat{v}}(X^{\epsilon,\hat{v}}(\tau_{D}^{\epsilon,\hat{v}})) \Bigr\} \notag \\
&= g_{k_0} (\bar{y}_0)  \label{Eq60} 
\end{align}
for $k \in \{1, 2, \ldots, n\}$, with $k_0 \in \{1, 2, \ldots, n\}$ and $\bar{y}_0 \in \partial D$. This completes the proof. 
\qed

Before concluding this section, it seems that appropriate to make a few more remarks concerning the implication of the main results. 

\begin{itemize}
\item Note that, from Proposition~\ref{P2}, the exit place on the boundary $\partial D$ for the component process $X^{\epsilon,\hat{v}}(t)$ starting from any $x \in D$ and random $k \in \{1, 2, \ldots, n\}$, satisfies
\begin{align*}
\mathbb{P}_{x,k}^{\epsilon,\hat{v}} \Bigl\{\bigl \vert X^{\epsilon,\hat{v}}(\tau_D^{\epsilon,\hat{v}}) - \bar{y}_0 \bigr \vert > \delta \Bigr\} \rightarrow 0 \quad \text{as} \quad \epsilon \rightarrow 0,
\end{align*}
that is, $X^{\epsilon,\hat{v}}(\tau_D^{\epsilon,\hat{v}}) \rightarrow \bar{y}_0$ in probability as $\epsilon \rightarrow 0$, for any $\delta > 0$, where $\bar{y}_0 \in \partial D$ is the unique minimum for the function $V$, and while the average dynamical system $\dot{x}(t) = \sum\nolimits_{k=1}^n \omega_k(x(t)) f_k \bigl(x(t), \hat{v}_k^{0}(x(t))\bigr)$ satisfies the large deviations condition (cf. Assumption~\ref{AS1}). Moreover, the type of distribution at the exit time $\tau_D^{\epsilon,\hat{v}}$ also satisfies
\begin{align*}
\mathbb{P}_{x,k}^{\epsilon,\hat{v}} \Bigl\{\zeta^{\epsilon}(\tau_D^{\epsilon,\hat{v}}) = k_0 \,\bigl \vert \tau_D^{\epsilon,\hat{v}}< \infty \Bigr\} = 1\rightarrow 0 \quad \text{as} \quad \epsilon \rightarrow 0,
\end{align*}
provided that there exists $k_0 \in \{1, 2, \ldots, n\}$ satisfying Assumption~\ref{AS2}. 
\item On the other hand, from Proposition~\ref{P1} (cf. Propositions~\ref{P3} and \ref{P4}), if one were to sit at the minimum of the function $V$, that is, the extremal point $\bar{y}_0 \in \partial D$ of the variational problem
\begin{align*}
V(y) = \inf \Bigl\{ I_{0T}(\varphi) \, \bigl\vert \, \varphi(0) = x_0 \in D, \,\, T > 0, \,\, \varphi(T) = y \,\,\, \text{for} \,\,\, y \in \partial D \Bigr\},
\end{align*}
where $V$ attained its unique minimum (i.e., $V(\bar{y}_0) < V(y)$, for any $y \in \partial D$), and wait for the process $(X^{\epsilon,\hat{v}}(t), r^{\epsilon,\hat{v}}(t))$ starting from any $x \in D$ and random $k \in \{1, 2, \ldots, n\}$ to arrive, then trace its history back along with all changes in the modes of the dynamics of the controlled diffusions. Furthermore, one would be likely to find that the process $(X^{\epsilon,\hat{v}}(t), r^{\epsilon,\hat{v}}(t))$ followed the most efficient path, as if it were guided. Note that the corresponding large deviations interpretation is simply that the dynamical system with even small random perturbations, will experiment with every conceivable paths, with probabilities that are extremely small. The higher the rate function $I_{0T}$, that is, the less efficient path, the smaller will be the probability and, hence, the less frequent attempts involving that path. Therefore, the first attempt to take the path to its minimum location will be the most efficient path.
\item Notice that, the asymptotic results for the exit place and the type of distribution at the exit time (or just before exiting D) are interpreted as giving an extremal path for the rate function $I_{0T}$ connecting $x_0 \in D$ and $\bar{y}_0 \in \partial D$, w.r.t. the admissible Markov controls $\hat{v}_{k}^{\epsilon} \in \mathcal{U}_k^{\rm ad}$, $k=1, 2, \ldots, n$.\footnote{Note that the last segment of the trajectory $X^{\epsilon,\hat{v}}(t)$ just before exiting $D$ belongs to a small neighborhood of this extremal path with probability close one as $\epsilon$ small enough.} Moreover, such admissible controls are also assumed rendering a non-empty closed invariant set in $\bar{D}$, when $\epsilon \rightarrow 0$, for the corresponding dynamical system $\dot{x}(t) = f_{k}\bigr(x(t), \hat{v}_{k}^{0}(x(t))\bigr)$, $k \in \{1,2, \ldots, n\} $ (see also Remark~\ref{R1} for the asymptotic behavior of the corresponding singular eigenvalue problem). Finally, from Proposition~\ref{P5}, one also observes that there is an intimate connection between the asymptotic results at the exit time and that of the limiting behavior of the Dirichlet problem for the corresponding system of elliptic PDEs with small vanishing parameter.
\end{itemize}


\end{document}